# SOME OPTIMALITY CONDITIONS OF SET-VALUED OPTIMIZATION PROBLEMS IN LOCALLY CONVEX TOPOLOGICAL VECTOR SPACES

Renying Zeng
Mathematics Department, Saskatchewan Polytechnic, Saskatoon, SK Canada S7L 4J7
Email: renying.zeng@saskpolytech.ca

**Abstract**: In this article, we work with set-valued optimization problems in locally convex topological vector spaces. We prove the equivalencies of some definitions of generalized convex maps introduced by Jeyakumar, Yang, Yang & Yang & Chen, as well as Zeng. And then, we discuss the conditions of weakly efficient solutions, proper efficient solutions, and optimal solutions of set-valued optimization problems.

## 1. Introduction and Preliminary

Leonhard Euler, in the eighteenth century, stated that "nothing at all takes place in the universe in which some rule of maximum and minimum does not appear." The statement may strike us extremely, yet it is undeniable that humankind's endeavors at least are usually associated with a quest for an optimum.

The development of optimization methods has been of interest in mathematics for centuries. The theory of mathematical optimization or mathematical programming is at the crossroads of many subjects. The subject grew from a realization that quantitative problems in manifestly different disciplines have important mathematical elements in common. Because of this commonality, many problems can be formulated and solved by using the unified set of ideas and methods that make up the field of optimization. The terms "minimum," "maximum," and "optimum" are in line with the mathematical tradition. Historically, linear programs were the focus in the optimization community, and initially, it was thought that the major divide was between linear and nonlinear optimization problems; later, people discovered that some nonlinear problems were much harder than others, and the "right" divide was between convex and nonconvex problems.

*During the years, many of our colleague mathematicians made contributions in generalized convexities and generalized convex optimizations, including well-known American mathematician Ky Fan* (see www.en.wikipedia.org) *as well as the international academic organization Working Group on Generalized Convexity (*WGGC*)-www.genconv.org.*

A topological vector space $Y$ with a convex cone $Y_+$ is said to be an ordered topological vector space, the partial order in $Y$ is given by:
$$y^1 \geq y^2, \text{ iff } y^1 - y^2 \in Y_+;$$

$$y^1 > y^2, \text{ iff } y^1 - y^2 \in \text{int } Y_+,$$
where $\text{int } Y_+$ is the interior of $Y_+$.

*In this article, we assume that $\text{int } Y_+ \neq \emptyset$.*

A subset $Y_+$ of $Y$ is said to be a cone if $\lambda y \in Y_+$ for all $y \in Y_+$ and all real scalars $\lambda \geq 0$. We denote by $0_Y$ the zero element in the topological vector space $Y$ and simply by $0$ if there is no confusion. A convex cone is one that satisfies $\lambda_1 y_1 + \lambda_2 y_2 \in Y_+$ for all $y_1, y_2 \in Y_+$ and all real scalar $\lambda_1, \lambda_2 \geq 0$. A pointed cone is one that satisfies $Y_+ \cap (-Y_+) = \{0\}$.

Suppose that $Y^*$ is the topological dual of $Y$, then
$$Y_+^* = \{\xi \in Y^* : \xi(y) \geq 0, \forall y \in Y_+\}$$
is said to be the dual cone of $Y_+$.

Let $X$ is a nonempty set and $D$ is a nonempty subset of $X$.

Suppose $f: X \to 2^Y$ is a set-valued function, where $X$ and $Y$ are topological vector spaces, $2^Y$ denotes the power set of $Y$. For any nonempty $D \subseteq X$, let
$$f(D) = \bigcup_{x \in D} f(x).$$

For $x \in X, \xi \in Y^*$, we set
$$\xi(f(x)) \geq 0, \text{ iff } \xi(y) \geq 0, \forall y \in f(x);$$
$$\xi(f(D)) \geq 0, \text{ iff } \xi(f(x)) \geq 0, \forall x \in D.$$

For any subset $A, B$ of $X$ and $u \in X$,
$$u + A = \{u + x : x \in A\},$$
$$A + B = \{x + y : x \in A, y \in B\}.$$

We denote by $R$ the set of all real numbers. For $A, B \subseteq R$, write
$$A \geq B, \text{ if } a \geq b, \forall a \in A, \forall b \in B.$$

A subset $M$ of $X$ is said to be convex, if $x_1, x_2 \in M$ and $0 < \alpha < 1$ implies $\alpha x_1 + (1-\alpha) x_2 \in M$; $M$ is said to be balanced if $x \in M$ and $|\alpha| \leq 1$ implies $\alpha x \in M$; $M$ is said to be absorbing if for any given neighborhood $N$ of $0$, there exists a positive scalar $\beta$ such that $\beta^{-1} M \subseteq U$, where $\beta^{-1} M = \{x \in X; x = \beta^{-1} v; v \in M\}$.

A topological vector space $X$ is called a locally convex topological space if any neighborhood $N$ of $0_X$ contains a convex, balanced, and absorbing open set.

*Normed linear spaces are locally convex topological spaces.*

## 2. The Equivalencies of Some Definitions of Generalized Convexities

*The well-known Chinese-American mathematician Ky Fan made fundamental contributions to operator & matrix theory, convex analysis & inequalities, linear & nonlinear programming, topology & fixed-point theory, and topological groups (cited from www.en.wikipedia.org). Fan defined the generalized convexity for vector-valued functions as follows* [1].

A function $f: D \subseteq X \to Y$ is said to be convexlike on $D$ if $\forall x_1, x_2 \in D, \forall \alpha \in (0,1)$ $\exists x_3 \in D$ such that

$$\alpha f(x_1) + (1-\alpha) f(x_2) \leq f(x_3).$$

For *set-valued maps*, the definition of "convexlike" can be modified as follows.

A function $f: X \to 2^Y$ is said to be convexlike on $D \subseteq X$ if $\forall x_1, x_2 \in D, \forall \alpha \in (0,1)$ $\exists x_3 \in D$ such that

$$\alpha f(x_1) + (1-\alpha) f(x_2) \subseteq f(x_3) + Y_+.$$

Jeyakumar, Yang, Yang et al., and Zeng [2-7] introduced the following Definitions 1-8. We proved in this paper that Definitions 1-4, and 5-8 are equivalent, respectively.

*I believe these equivalences are very interesting results in recent research of generalized convexities. It even surprised myself since the proofs were so straight forward for locally convex topological vector spaces by using the concepts of balanced, absorbent, and convex sets..*

**Definition 1** A set-valued map $f : X \to 2^Y$ is called to be subconvexlike on $D \subseteq X$ if
$$\exists \theta \in \text{int } Y_+, \forall x^1, x^2 \in D, \forall \alpha \in (0,1), \forall \varepsilon > 0, \exists x^3 \in D, s.t.$$
$$\varepsilon \theta + \alpha f(x^1) + (1-\alpha) f(x^2) \subseteq f(x^3) + Y_+.$$

**Definition 2** A set-valued map $f : X \to 2^Y$ is called to be sub-convexlike on $D \subseteq X$ if
$$\exists u \in U, \forall x^1, x^2 \in D, \forall \alpha \in (0,1), \forall \varepsilon > 0, \exists x^3 \in D, s.t.$$
$$\varepsilon u + \alpha f(x^1) + (1-\alpha) f(x^2) \subseteq f(x^3) + Y_+,$$
where $U$ is the set of all bounded functions from $X$ to $2^Y$.

**Definition 3** A set-valued map $f : X \to 2^Y$ is called to be strictly subconvexlike on $D \subseteq X$ if
$$\exists \theta \in \text{int } Y_+, \forall x^1, x^2 \in D, \forall \alpha \in (0,1), \forall \varepsilon > 0, \exists x^3 \in D, s.t.$$
$$\varepsilon \theta + \alpha f(x^1) + (1-\alpha) f(x^2) \subseteq f(x^3) + \text{int } Y_+.$$

**Definition 4** A set-valued map $f : X \to 2^Y$ is called to be strictly sub-convexlike on $D \subseteq X$ if
$$\exists u \in U, \forall x^1, x^2 \in D, \forall \alpha \in (0,1), \forall \varepsilon > 0, \exists x^3 \in D, s.t.$$
$$\varepsilon u + \alpha f(x^1) + (1-\alpha)f(x^2) \subseteq f(x^3) + \text{int } Y_+.$$

**Theorem 1** Definitions 1, 2, 3, and 4 are equivalent.

**Proof.** a) Definition 1 $\Rightarrow$ Definition 2. We take $u \equiv \{\theta\} : D \to 2^Y$ a bounded set-valued map. So Definition 1 implies Definition 2.

b) Definition 2 $\Rightarrow$ Definition 3.

Because $f$ is subconvexlike,
$$\exists u \in U, \forall x^1, x^2 \in D, \forall \alpha \in (0,1), \forall \delta > 0, \exists x^3 \in D, s.t.$$
$$\delta u + \alpha f(x^1) + (1-\alpha)f(x^2) \subseteq f(x^3) + Y_+,$$

It is known that, the interior $\text{int } Y_+$ of a convex cone $Y_+$ is also a convex cone. Therefore, for any given $\theta \in \text{int } Y_+$, and any given $\varepsilon > 0$ we have $\varepsilon\theta \in \text{int } Y_+$. It follows also that there exists a neighborhood $N$ of the origin $0_Y$ of $Y$ such that $\varepsilon\theta + N \subseteq \text{int } Y_+$. Since $Y$ is a locally convex space, without loss of generality, we may assume that $N$ is balanced and absorbent. We can choose a positive number $\delta$ which is sufficiently small such that $-\delta u \in N$. Hence
$$\varepsilon\theta - \delta u \in \text{int } Y_+.$$

Therefore,
$$\alpha f(x^1) + (1-\alpha)f(x^2) \subseteq f(x^3) - \delta u + Y_+.$$

Hence
$$\varepsilon\theta + \alpha f(x^1) + (1-\alpha)f(x^2)$$
$$\subseteq \varepsilon\theta + f(x^3) - \delta u + Y_+$$
$$\subseteq f(x^3) + \text{int } Y_+ + Y_+$$
$$\subseteq f(x^3) + \text{int } Y_+.$$

c) Definition 3 $\Rightarrow$ Definition 4 is obvious.

d) Definition 4 $\Rightarrow$ Definition 1. The proof is similar to b). □

**Definition 5** A set-valued map $f : X \to 2^Y$ is called to be presubconvexlike on $D \subseteq X$ if
$$\exists \theta \in \text{int } Y_+, \forall x^1, x^2 \in D, \forall \alpha \in (0,1), \forall \varepsilon > 0, \exists x^3 \in D, \exists \tau > 0, s.t.$$
$$\varepsilon\theta + \alpha f(x^1) + (1-\alpha)f(x^2) \subseteq \tau f(x^3) + Y_+.$$

**Definition 6** A set-valued map $f : X \to 2^Y$ is called to be presub-convexlike on $D \subseteq X$ if
$$\exists u \in U, \forall x^1, x^2 \in D, \forall \alpha \in (0,1), \forall \varepsilon > 0, \exists x^3 \in D, \exists \tau > 0, s.t.$$
$$\varepsilon u + \alpha f(x^1) + (1-\alpha)f(x^2) \subseteq \tau f(x^3) + Y_+,$$
where $U$ is the set of all bounded functions from $X$ to $2^Y$.

**Definition 7** A set-valued map $f : X \to 2^Y$ is called to be strictly presubconvexlike on $D \subseteq X$ if
$$\exists \theta \in \text{int } Y_+, \forall x^1, x^2 \in D, \forall \alpha \in (0,1), \forall \varepsilon > 0, \exists x^3 \in D, \exists \tau > 0, s.t.$$
$$\varepsilon \theta + \alpha f(x^1) + (1-\alpha)f(x^2) \subseteq \tau f(x^3) + \text{int } Y_+.$$

**Definition 8** A set-valued map $f : X \to 2^Y$ is called to be strictly presub-convexlike on $D \subseteq X$ if
$$\exists u \in U, \forall x^1, x^2 \in D, \forall \alpha \in (0,1), \forall \varepsilon > 0, \exists x^3 \in D, \exists \tau > 0, s.t.$$
$$\varepsilon u + \alpha f(x^1) + (1-\alpha)f(x^2) \subseteq \tau f(x^3) + \text{int } Y_+.$$

**Theorem 2** Definitions 5, 6, 7, and 8 are equivalent.

**Proof.** e) Definition 5 $\Rightarrow$ Definition 6. Note that $u \equiv \{\theta\}: D \to 2^Y$ is a bounded set-valued map. So, Definition 5 implies Definition 6.

f) Definition 6 $\Rightarrow$ Definition 7.

Because $f$ is presubconvexlike,
$$\exists u \in U, \forall x^1, x^2 \in D, \forall \alpha \in (0,1), \forall \delta > 0, \exists x^3 \in D, \exists \tau > 0, s.t.$$
$$\delta u + \alpha f(x^1) + (1-\alpha)f(x^2) \subseteq \tau f(x^3) + Y_+.$$

Given $\theta \in \text{int } Y_+$. For any given $\varepsilon > 0$ we have $\varepsilon \theta \in \text{int } Y_+$. It follows that there exists a neighborhood $N$ of the origin $0_Y$ of $Y$ such that $\varepsilon \theta + N \subseteq \text{int } Y_+$. Since Y is a locally convex space, we may assume that $N$ is balanced and absorbent. So, we may choose positive number $\delta$ which is sufficiently small such that $-\delta u \in N$. Hence
$$\varepsilon \theta - \delta u \in \text{int } Y_+.$$
Therefore, $\delta u + \alpha f(x^1) + (1-\alpha)f(x^2) \subseteq \tau f(x^3) + Y_+$ given that
$$\alpha f(x^1) + (1-\alpha)f(x^2) \subseteq \tau f(x^3) - \delta u + Y_+.$$
Consequently,
$$\varepsilon \theta + \alpha f(x^1) + (1-\alpha)f(x^2)$$
$$\subseteq \varepsilon \theta + \tau f(x^3) - \delta u + Y_+$$
$$\subseteq \tau f(x^3) + \text{int } Y_+ + Y_+$$
$$\subseteq \tau f(x^3) + \text{int } Y_+$$

g) Definition 7 ⇒ Definition 8 is obvious.

h) Definition 8 ⇒ Definition 5. The proof is similar to f ). □

## 3. Weakly Efficient Solutions

Consider the following optimization problem with set-valued maps:

$$(VP) \quad \begin{array}{l} Y_+ - \min f(x), \\ s.t., g(x) \cap (-Z_+) \neq \emptyset, 0_W \in h(x), \\ x \in X. \end{array}$$

Where $X$, $Y$, $Z$ and $W$ are all locally convex topological vector spaces.

From now on, let $D$ be the feasible set of (VP), i.e.,

$$D = \{x \in X : g(x) \cap (-Y_+) \neq \emptyset, 0_W \in h(x)\}.$$

**Definition 9** $\bar{x} \in D$ is said to be a weakly efficient solution of (VP), if $\exists \bar{y} \in f(\bar{x})$ such that $\forall x \in D$, there is no $y \in f(x)$ satisfying $(\bar{y} - y) \in \operatorname{int} Y_+$.

It is known that $\bar{x} \in D$ is a weakly efficient solution of (VP) if and only if $\exists \bar{y} \in f(\bar{x})$ such that $(\bar{y} - f(D)) \cap \operatorname{int} Y_+ = \emptyset$.

**Definition 10** The problem (VP) is said to satisfy the Slater constraint qualification (SC) if $\forall (\eta, \varsigma) \in (Z_+^* \times W^*) \setminus \{(0_{Z^*}, 0_{W^*})\}, \exists x \in D$ and $\exists$ a negative real number $t \in \eta(g(x)) \cap \varsigma(h(x))$.

Let

$$P\min[A, Y_+] = \{y \in A : (y - A) \cap \operatorname{int} Y_+ = \emptyset\},$$

$$P\max[A, Y_+] = \{y \in A : (A - y) \cap \operatorname{int} Y_+ = \emptyset\}.$$

**Definition 11** A triple $(\bar{x}, \bar{S}, \bar{T}) \in X \times B^+(Z, Y) \times B(W, Y)$ is said to be a weak saddle point of $L$ if

$$L(\bar{x}, \bar{S}, \bar{T}) \cap P\min[L(X, \bar{S}, \bar{T}), Y_+]$$
$$\cap P\max[L(\bar{x}, B^+(Z, Y), B(W, Y)), Y_+]$$
$$\neq \emptyset.$$

Where

$$L(\bar{x},\bar{S},\bar{T}) = f(\bar{x}) + \bar{S}(g(\bar{x})) + \bar{T}(h(\bar{x})).$$

*Zeng [8] introduced a definition of prenearaffinelikeness for vector-valued maps by using of "affine cones". The following Definition 12 is a definition of generalized affinelike set-valued maps, but use of "linear subspaces".*

**Definition 12** Suppose $E \subseteq Y$ is a linear subspace of Y. A set-valued map $f: X \to 2^Y$ is said to be $E$-generalized affinelike on $D$ if $\forall x_1, x_2 \in D$, $\forall \alpha \in R$, $\exists v \in E$, $\exists x_3 \in D, \exists \tau > 0$ such that

$$v + \alpha f(x_1) + (1-\alpha)f(x_2) \subseteq \tau f(x_3).$$

**Lemma 1** Suppose that $f(x), g(x)$, and $h(x)$ are all set-valued maps, $E$ is a linear subspace of $W$. If $f(x), g(x)$ are presubconvexlike on $D$, $h(x)$ is $E$-generalized affinelike on $D$ and $int[h(D)] \neq \emptyset$, and if (i) and (ii) denote the systems

(i) $\exists x \in D, s.t., f(x) \cap (-int\, Y_+) \neq \emptyset$, $g(x) \cap (-Z_+) \neq \emptyset, 0_W \in h(x)$;

(ii) $\exists (\xi, \eta, \varsigma) \in (Y_+^* \times Z_+^* \times W^*) \setminus \{(0_Y, 0_Z, 0_W)\}$ such that

$$\xi(f(x)) + \eta(g(x)) + \varsigma(h(x)) \geq 0, \ \forall x \in D.$$

If (i) has no solution then (ii) has solutions.

Moreover if (ii) has a solution $(\xi, \eta, \varsigma)$ with $\xi \neq 0_{Y^*}$ then (i) has no solutions.

**Proof.** Similar to the proof of Theorem 4.1 in [8].

**Theorem 3** Suppose that $f(x), g(x)$, and $h(x)$ are all set-valued maps, $E$ is a linear subspace of $W$. Let $\bar{x} \in D$ be weakly efficient solution of (VP). If $f(x) - f(\bar{x}), g(x)$ are presubconvexlike on $D$, $h$ $(x)$ is $E$-generalized affinelike on $D$ and $int[h(D)] \neq \emptyset$, and if (VP) satisfies the Slater constrained qualification (SC), then $\exists (\bar{S}, \bar{T}) \in B^+(Z, Y) \times B(W, Y)$ such that $(\bar{x}, \bar{S}, \bar{T}) \in X \times B^+(Z, Y) \times B(W, Y)$ is a weak saddle point of $L$ and $O \in \bar{S}(g(\bar{x}))$, $O \in \bar{T}(h(\bar{x}))$.

**Proof**. Suppose $\bar{x} \in D$ is a weakly efficient solution of (VP), then $\exists \bar{y} \in f(\bar{x})$ for which there is not any $x \in D$ such that $f(x) - \bar{y} \in -int\, Y_+$. So, there is not any $x \in D$ such that $f(x) - \bar{y} \in -int\, Y_+, g(x) \in -Z_+, O \in h(x)$. Then, by Lemma 1, $\exists (\xi, \eta, \varsigma) \in Y_+^* \times Z_+^* \times W^* \setminus \{O\}$ such that

$$\xi(f(x)-\bar{y})+\eta(g(x))+\varsigma(h(x))\geq 0, \forall x\in D.$$

Then, we conclude that

$$O\in\bar{S}(g(\bar{x})), O\in\bar{T}(h(\bar{x})),$$

as well as $(\bar{y}-[f(D)+S(g(D))+T(h(D))])\cap \text{int } Y_+ =\emptyset$. Hence

$$\bar{y}\in P\min[L(X,\bar{S},\bar{T}),Y_+].$$

On the other hand, since $O\in\bar{S}(g(\bar{x}))$, $\exists \bar{z}\in g(\bar{x})$ we get $\bar{S}(\bar{z})=O$. This and $(f(\bar{x})-\bar{y})\cap \text{int } Y_+ =\emptyset$ together deduce that

$$(f(\bar{x})-\bar{y}-S(\bar{z}))\cap \text{int } Y_+ =\emptyset.$$

Therefore, we conclude that $(\bar{x},\bar{S},\bar{T})\in X\times B^+(Z,Y)\times B(W,Y)$ is a weak saddle point of $L$. □

**Theorem 4** Suppose that $f(x), g(x)$, and $h(x)$ are all set-valued maps. A triple $(\bar{x},\bar{S},\bar{T})\in X\times B^+(Z,Y)\times B^+(W,Y)$ is a weak saddle point of $L$ if and only if $\exists \bar{y}\in f(\bar{x}), \bar{z}\in g(\bar{x}), \bar{w}\in h(\bar{x})$, such that

(i) $\bar{y}\in P\min[\bigcup_{x\in X} L(x,\bar{S},\bar{T}),Y_+]\cap P\max[f(\bar{x}),Y_+]$,

(ii) $-g(\bar{x})\subseteq Y_+$,

(iii) $0_Y \in \bar{T}(\bar{z})$.

**Proof**. We only need to prove the necessity. Since $(\bar{x},\bar{S},\bar{T})\in X\times B^+(Z,Y)\times B^+(W,Y)$ is a weak saddle point of $L$, $\exists \bar{y}\in f(\bar{x}), \bar{z}\in g(\bar{x}), \bar{w}\in h(\bar{x})$ such that

$$\bar{y}+\bar{S}(\bar{z})+\bar{T}(\bar{w})\in P\min[\bigcup_{x\in X}[L(x,\bar{S},\bar{T}),Y_+],$$

and

$$\bar{y}+\bar{S}(\bar{z})+\bar{T}(\bar{w})\in P\max[\bigcup_{(S,T)\in B^+(Z,Y)\times B^+(W,Y)} L(\bar{x},S,T),Y_+].$$

then we deduce that

$$(f(\bar{x})+S(g(\bar{x}))+T(h(\bar{x}))-[\bar{y}+\bar{S}(\bar{z})+\bar{T}(\bar{w})])\cap \text{int } Y_+ =\emptyset,$$
$$\forall (S,T)\in B^+(Z,Y)\times B^+(W,Y).$$

Then, $\exists \bar{y}\in f(\bar{x})$, for $\forall x\in D$, there exists no $y\in f(x)$ such that $y-\bar{y}\in -\text{int } Y_+$. By

$$\xi[f(x)+S(g(x))+T(h(x))]$$
$$\geq \xi(\bar{y}), \forall x\in D.$$

We have

$$\bar{y}\in f(\bar{x})\subseteq f(\bar{x})+S(g(\bar{x}))+T(h(\bar{x})).$$

Therefore
$$\bar{y} \in P\min[\bigcup_{x \in X} L(x,\bar{S},\bar{T}),Y_+] \cap P\max[f(\bar{x}),Y_+]. \square$$

## 4. Properly Efficient Solutions

**Definition 13** Given $\bar{\xi} \in Y_+^* \setminus \{0_{Y^*}\}$. The real-valued Lagrangian mapping of (VP) $l_{\bar{\xi}} : X \times Z_+^* \times W^* \to R$ is defined by
$$l_{\bar{\xi}}(x,\eta,\varsigma) = \bar{\xi}(f(x)) + \eta(g(x)) + \varsigma(h(x)).$$

**Definition 14** Given $\bar{\xi} \in Y_+^* \setminus \{0_{Y^*}\}$. A triple $(\bar{x},\bar{\eta},\bar{\varsigma})$ is said to be a saddle point of the Lagrangian mapping $l_{\bar{\xi}}$, if
$$l_{\bar{\xi}}(\bar{x},\eta,\varsigma) \leq l_{\bar{\xi}}(\bar{x},\bar{\eta},\bar{\varsigma}) \leq l_{\bar{\xi}}(x,\bar{\eta},\bar{\varsigma}),$$
for $\forall x \in D, \forall(\eta,\varsigma) \in Z_+^* \times W^*$.

**Definition 15** A vector $\bar{x} \in D$ is said to be a properly efficient solution of (VP) if $\exists \xi \in Y_+^* \setminus \{0_{Y^*}\}$ such that
$$\xi(f(x)) \geq \xi(f(\bar{x})).$$

**Theorem 5** Suppose that $f(x), g(x)$, and $h(x)$ are all set-valued maps, $E$ is a linear subspace of $W$. Let $\bar{x} \in D$. If $\exists \bar{\xi} \in Y_+^* \setminus \{0_{Y^*}\}$ and $\exists (\bar{\eta},\bar{\varsigma}) \in Z_+^* \times W^*$ such that $(\bar{x},\bar{\eta},\bar{\varsigma})$ is a saddle point of the Lagrangian mapping $l_{\bar{\xi}}$, then $\bar{x} \in D$ is a properly efficient solution of (VP) and $\bar{\eta}(g(\bar{x})) = \{0\}$, $\bar{\varsigma}(h(\bar{x})) = \{0\}$.

**Proof.** Suppose $(\bar{x},\bar{\eta},\bar{\varsigma})$ is a saddle point of the Lagrangian mapping $l_{\bar{\xi}}$, then
$$\eta(g(\bar{x})) + \varsigma(h(\bar{x}))) \leq \bar{\eta}(g(\bar{x})) + \bar{\varsigma}(h(\bar{x})), \forall (\eta,\varsigma) \in Z_+^* \times W^*.$$
Take $\eta = \bar{\eta}$, or $\varsigma = \bar{\varsigma}$ we have
$$\eta(g(\bar{x})) \leq \bar{\eta}(g(\bar{x})), \forall \eta \in Z_+^*,$$
$$\varsigma(h(\bar{x}))) \leq \bar{\varsigma}(h(\bar{x})), \forall \varsigma \in W^*.$$
Taking $\eta = 0_{Z^*}$ we get $\bar{\eta}(g(\bar{x})) \geq 0$, but taking $\eta = 2\bar{\eta}$ we have $\bar{\eta}(g(\bar{x})) \leq 0$. Hence
$$\bar{\eta}(g(\bar{x})) = \{0\}.$$
Similarly,
$$\bar{\varsigma}(h(\bar{x})) = \{0\}.$$
Noting that $g(x) \cap (-Z_+) \neq \emptyset$ and $0 \in \bar{\varsigma}(h(x))$ we have
$$\bar{\xi}(f(\bar{x})) \leq \bar{\xi}(f(x)), \forall x \in D.$$
Therefore, $\bar{x}$ is a properly efficient solution of (VP). $\square$

**Theorem 6** Suppose that $f(x), g(x)$, and $h(x)$ are all set-valued maps, $E$ is a linear subspace of $W$. Let $\bar{x} \in D$ be properly efficient solution of (VP). If $f(x) - f(\bar{x}), g(x)$ are presubconvexlike on $D$, $h(x)$ is $E$-generalized affinelike on $D$ and $int[h(D)] \neq \emptyset$, and if (VP) satisfies the Slater constrained qualification (SC), then $\exists \bar{\xi} \in Y_+^* \setminus \{0_{Y^*}\}$, and $\exists (\bar{\eta}, \bar{\varsigma}) \in Z_+^* \times W^*$ such that $(\bar{x}, \bar{\eta}, \bar{\varsigma})$ is a saddle point of the Lagrangian mapping $l_{\bar{\xi}}$ and $\bar{\eta}(g(\bar{x})) = \{0\}$, $\bar{\eta}(g(\bar{x})) = \{0\}$.

## 5. Optimal Solutions

**Definition 16** The Vector Lagrangian map $L: X \times B^+(Z, Y) \times B^+(W, Y) \to 2^Y$ for (VP) is defined to be the set-valued map
$$L(x, S, T) = f(x) + S(g(x)) + T(h(x)).$$

Given $(S, T) \in B^+(Z, Y) \times B^+(W, Y)$. We consider the following unconstrained vector minimization problem induced by (VP):

(VPST1) $\quad Y_+ - \min L(x, S, T),$
$\quad\quad\quad\quad\quad s.t., x \in D.$

**Theorem 7** Suppose that $f(x), g(x)$, and $h(x)$ are all set-valued maps, $E$ is a linear subspace of $W$. Let $\bar{x} \in D$. If $f(x) - f(\bar{x}), g(x)$ are presubconvexlike on $D$, $h(x)$ is $E$-generalized affinelike on $D$ and $int[h(D)] \neq \emptyset$, and if (VP) satisfies the Slater constrained qualification (SC), then $\bar{x}$ is a weakly efficient solution of (VP) if and only if $\exists (S, T) \in B^+(Z, Y) \times B^+(W, Y)$ such that $\bar{x}$ is an optimal solution of (VPST1) and $0_Y \in S(g(\bar{x}))$, $0_Y \in T(h(\bar{x}))$.

**Proof.** Suppose that $\bar{x} \in D$ is a weakly efficient solution of (VP), then $\exists \bar{y} \in f(\bar{x})$ there is not any $x \in D$ such that $f(x) - \bar{y} \in -int Y_+$. That is to say, there is not any $x \in X$ such that
$$f(x) - \bar{y} \in -int Y_+, g(x) \in -Z_+, 0_W \in h(x).$$
Therefore, $\exists (\xi, \eta, \varsigma) \in Y_+^* \setminus \{0_{Y^*}\} \times Z_+^* \times W^*$ such that
$$\xi(f(x) - \bar{y}) + \eta(g(x)) + \varsigma(h(x)) \geq 0, \forall x \in D.$$
Taking $x = \bar{x}$ in the above we obtain
$$\xi(f(x) - \bar{y}) \geq 0, \forall x \in D.$$
Since $\xi \neq 0_{Y^*}$ we may take $y_0 \in Y_+ \setminus \{0_Y\}$ such that
$$\xi(y_0) = 1.$$
Define the operator $S: Z \to Y$ and $T: W \to Y$ by

$$S(z) = \eta(z)y_0, T(w) = \varsigma(w)y_0.$$

It is easy to see that
$$S \in B^+(Z,Y), S(Z_+) = \eta(Z_+)y_0 \subseteq Y_+,$$
$$T \in B^+(W,Y), T(Z_+) = \varsigma(Z_+)y_0 \subseteq Y_+,$$

and
$$S(g(\bar{x})) = \eta(g(\bar{x}))y_0 \in 0 \cdot Y_+ = 0_Y.$$

Since $\bar{x} \in D$, $0_W \in h(\bar{x})$. Hence
$$0_Y \in T(h(\bar{x})).$$

Therefore
$$S \in B^+(Z,Y), \bar{y} \in f(\bar{x}) \subseteq f(\bar{x}) + S(g(\bar{x})) + T(h(\bar{x})).$$

And so
$$\xi[f(x) + S(g(x)) + T(h(x))]$$
$$= \xi(f(x)) + \eta((g(x))\xi(y_0) + \varsigma(h(x))\xi(y_0)$$
$$= \xi(f(x)) + \eta(g(x)) + \varsigma(h(x))$$
$$\geq \xi(\bar{y}), \forall x \in D.$$

i.e.,
$$\xi(f(x) - \bar{y}) + (\xi \circ S)(g(x)) + (\xi \circ T)(h(x))$$
$$\geq 0, \forall x \in D.$$

Because the compound operators $\xi \circ S \in Z_+^*, \xi \circ T \in W_+$, we conclude that, $\bar{x} \in D$ is an optimal solution of (VPST1).

On the other hand, assume that $\exists (S,T) \in B^+(Z,Y) \times B^+(W,Y)$ such that $0_Y \in S(g(\bar{x}))$, $0_Y \in T(h(\bar{x}))$ and $\bar{x} \in D$ is an optimal solution of (VPST1).

If $\bar{x} \in D$ is not weakly efficient solution of (VP), by the definition of weakly efficient solution, $\forall y \in f(\bar{x})$ we have
$$(y - f(D)) \cap \text{int } Y_+ \neq \emptyset.$$

Let $y - y_0 \in ((y - f(D)) \cap \text{int } Y_+$, where $y_0 \in f(D)$.

Since $0_Y \in S(g(\bar{x})), 0_Y \in T(h(\bar{x}))$, we get
$$v = y + 0_Z + 0_W = y \in f(\bar{x}) + S(g(\bar{x})) + T(h(\bar{x})),$$

and
$$v_0 = y_0 + 0_Z + 0_W = y_0$$
$$\in f(\bar{x}) + S(g(\bar{x})) + T(h(\bar{x}))$$
$$\subseteq f(D) + S(g(D)) + T(h(D)).$$

Therefore
$$[v - (f(D) + S(g(D)) + T(h(D)))] \cap \text{int } Y_+ \neq \emptyset,$$

which contradicts to the assumption that $\bar{x} \in D$ is an optimal solution of (VPST1). □

*Our results in Sections 3, 4 and 5 modified the corresponding results in [9-12], and our ideas and methods may be used to extend the results in [13, 14].*